\documentclass[12pt]{article}
\usepackage[a4paper,left=0.8in,right=0.8in,top=1in,bottom=1in]{geometry}

\usepackage[T1]{fontenc}
\usepackage[utf8]{inputenc}
\usepackage{lmodern}

\usepackage{amsmath,amssymb,amsthm}
\usepackage{graphicx}
\usepackage{tikz}
\usepackage{hyperref}
\usepackage{anyfontsize}
\usepackage{colortbl}
\usepackage{xcolor}
\usepackage{steinmetz}
\usepackage{setspace}
\graphicspath{{./images/}}

\title{On $t$-Edge-Balanced Graphs}
\author{Yeow Meng Chee\\
Singapore University of Technology and Design\\
8 Somapah Road, Singapore 487372\\
\texttt{ymchee@sutd.edu.sg}}
\date{April 2026}

\newtheorem{lemma}{Lemma}
\newtheorem{corollary}{Corollary}
\newtheorem{conjecture}{Conjecture}

\begin{document}

\maketitle

\begin{abstract}
A graph $G$ on $n$ vertices with $k$ edges is $t$-edge-balanced if every graph on $n$ vertices with $t$ edges is contained in exactly the same number of subgraphs of $K_n$ isomorphic to $G$. Despite the existence of infinite families of $2$-edge-balanced graphs, no $t$-edge-balanced graphs were known for $t \ge 3$. This paper resolves the existence question for $t \ge 3$ in two directions. For $t = 3$, we derive necessary arithmetic conditions on the parameters $(n,k)$ and use a simulated annealing search to find the first known examples of $3$-edge-balanced graphs. For $t \ge 4$, we prove that no nontrivial $t$-edge-balanced graphs exist.
\end{abstract}


\bigskip
\noindent\textbf{Funding:}
Research supported by SUTD Grant SKI 2021\_07\_04.

\section{Introduction}\label{sec1}

We work throughout with graphs on a fixed vertex set $[n]=\{1,\dots,n\}$. For an integer $k$, let $\mathcal{G}(n,k)$ denote the set of all graphs on $[n]$ with exactly $k$ edges. We consider the action of the symmetric group $S_n$ on $[n]$, which induces a natural action on the edge set of $K_n$, and hence on graphs with vertex set $[n]$. We follow standard graph-theoretic notation as in Bollob\'{a}s \cite{Bollobas:1998}.

Let $G \in \mathcal{G}(n,k)$, and let $0 < t < k$. We say that $G$ is \emph{$t$-edge-balanced} if there exists a positive integer $\lambda$, called the \emph{index} of $G$, such that every $H \in \mathcal{G}(n,t)$ is contained in exactly $\lambda$ subgraphs of $K_n$ that are isomorphic to $G$. Equivalently, every $H \in \mathcal{G}(n,t)$ extends to copies of $G$ with the same multiplicity. This notion is closely connected to graphical designs: $G$ is $t$-edge-balanced precisely when the $S_n$-orbit of $G$ forms the block set of a graphical $t$-design with point set $E(K_n)$ (see, for example, Chee and Kreher \cite{CheeKreher:2007}). Indeed, the notion of $2$-edge-balanced graphs was introduced by Alltop as a construction technique for $2$-designs \cite{Alltop:1966}.

The results below follow easily from the equivalence between $t$-edge-balanced graphs and single-orbit graphical $t$-designs.

\begin{lemma}\label{derived}
If $G$ is $t$-edge-balanced, then $G$ is also $s$-edge-balanced for every $s$, $0\le s\le t$.
\end{lemma}

\begin{lemma}
If $G\in\mathcal{G}(n,k)$ is $t$-edge-balanced, then its complement $K_n-G\in\mathcal{G}(n,\binom{n}{2}-k)$ is also $t$-edge-balanced.
\end{lemma}

Hence, in studying existence, we may restrict attention to graphs with $k \le \frac{1}{2}\binom{n}{2}$. We are also interested only in the nondegenerate range $2 \le t < k$. Indeed, the case $t=1$ is vacuous, since there is only one isomorphism class of graphs with one edge. The $t$-edge-balanced graphs satisfying $2 \le t < k\le \frac{1}{2}\binom{n}{2}$ are called \emph{nontrivial}.

The first nontrivial case is $t=2$. Alltop \cite{Alltop:1966} showed that the graph $G\in\mathcal{G}(2k-3,k)$ consisting of a single cycle of length $k$ and $k-3$ isolated vertices, is 2-edge-balanced for all $k\geq 4$. Further infinite families of $2$-edge-balanced graphs were obtained by Caliskan and Chee \cite{CaliskanChee:2014,Caliskan:2016}.

For $t\ge 3$, no $t$-edge-balanced graphs are known. In view of Lemma \ref{derived}, one would have hoped that among the known infinite families of 2-edge-balanced graphs, some turned out to be also 3-edge-balanced. However, a 3-edge-balanced graph must contain at least one copy of every 3-edge isomorphism class, in particular a triangle $K_3$, and all of the known infinite families of 2-edge-balanced graphs are triangle-free. 

This paper gives the first examples of $3$-edge-balanced graphs, establishing their existence. These examples include the 10 smallest $3$-edge-balanced graphs that can exist. For $t\geq 4$, we show that no $t$-edge-balanced graphs exist.

\section{Preliminaries}\label{sec2}

Recall that the symmetric group $S_n$ acts on $[n]$, which induces an action on the edge set of $K_n$, and hence on graphs with vertex set $[n]$. The orbits of this action are precisely the isomorphism classes of graphs on $[n]$. In particular, for any $t$, the action of $S_n$ partitions $\mathcal{G}(n,t)$ into isomorphism classes of $t$-edge graphs.

For graphs $G$ and $H$ on $[n]$, let $n_{G:H}$ denote the number of subgraphs of $G$ isomorphic to $H$, and let $\lambda_{G:H}$ denote the number of subgraphs of $K_n$ isomorphic to $G$ that contain $H$ as a subgraph. The key observation is that $\lambda_{G:H}$ can be computed directly from $n_{G:H}$ and the automorphism groups of $H$ and $G$.

\begin{lemma}[Alltop \cite{Alltop:1966}]\label{alltop}
Let $G$ and $H$ be graphs on $[n]$. Then
\[
\lambda_{G:H} = n_{G:H} \frac{|\mathrm{Aut}(H)|}{|\mathrm{Aut}(G)|}.
\]
\end{lemma}

The utility of Lemma~\ref{alltop} is as follows. To check whether $G$ is $t$-edge-balanced, one needs to verify that every $H \in \mathcal{G}(n,t)$ is contained in the same number of copies of $G$ in $K_n$. Since $\lambda_{G:H}$ depends on $H$ only through $n_{G:H}$ and $|\mathrm{Aut}(H)|$, this reduces to a condition that can be checked by examining the structure of $G$ alone, without enumerating copies of $G$ in $K_n$ directly.

More precisely, with $H_1, \ldots, H_r$ the isomorphism class representatives of $\mathcal{G}(n,t)$, $G \in \mathcal{G}(n,k)$ is $t$-edge-balanced if and only if $\lambda_{G:H_i}$ takes the same value for all $i \in \{1, \ldots, r\}$. By Lemma~\ref{alltop}, this is equivalent to
\begin{equation}\label{balanced-condition}
n_{G:H_i} \cdot |\mathrm{Aut}(H_i)| = n_{G:H_j} \cdot |\mathrm{Aut}(H_j)|,
\end{equation}
for all $i, j \in \{1, \ldots, r\}$. Thus, to construct a $t$-edge-balanced graph, it suffices to find a graph $G$ whose subgraph counts $n_{G:H_i}$ satisfy \eqref{balanced-condition}. When this holds, the common value of $\lambda_{G:H_i}$ is the index $\lambda$ of $G$, and it can be computed from Lemma~\ref{alltop} once the structure of $G$ is known.

\section{Necessary Conditions for 3-Edge-Balanced Graphs}\label{sec3}

Let $N=\binom{n}{2}$. As noted in Section~\ref{sec1}, we may restrict attention to $4\le k\le N/2$. In this section we derive two necessary arithmetic conditions on the parameter pair $(n,k)$ for a $3$-edge-balanced graph in $\mathcal{G}(n,k)$ to exist, and identify the smallest admissible parameter sets.

By Lemma~\ref{derived}, any $3$-edge-balanced graph is also $2$-edge-balanced, so we begin by deriving the constraints imposed by $2$-edge-balance. The two isomorphism classes of graphs in $\mathcal{G}(n,2)$ are listed in Table~\ref{tab:2edge}. Note that $2K_2$ requires $n\ge 4$; for $n<4$ there is only one isomorphism class and every graph is trivially 2-edge-balanced.

\begin{table}[h]
\centering
\begin{tabular}{llll}
\hline
Graph $H$ & Figure & $\mathrm{Aut}(H)$ & $|\mathrm{Aut}(H)|$ \\
\hline
$P_3$ & \begin{tikzpicture}[baseline=-0.5ex, vertex/.style={circle, fill=black, inner sep=1.5pt}]
\node[vertex] (1) at (0,0) {};
\node[vertex] (2) at (0.5,0) {};
\node[vertex] (3) at (1,0) {};
\draw (1)--(2)--(3);
\end{tikzpicture} & $S_2 \times S_{n-3}$ & $2(n-3)!$ \\
$2K_2$ & \begin{tikzpicture}[baseline=-0.5ex, vertex/.style={circle, fill=black, inner sep=1.5pt}]
\node[vertex] (1) at (0,0) {};
\node[vertex] (2) at (0.5,0) {};
\node[vertex] (3) at (1,0) {};
\node[vertex] (4) at (1.5,0) {};
\draw (1)--(2);
\draw (3)--(4);
\end{tikzpicture} & $(S_2 \wr S_2) \times S_{n-4}$ & $8(n-4)!$ \\
\hline
\end{tabular}
\caption{The two isomorphism classes of graphs in $\mathcal{G}(n,2)$ and their automorphism groups.}
\label{tab:2edge}
\end{table}

Let $a$ and $b$ denote the number of subgraphs of $G$ isomorphic to $P_3$ and $2K_2$, respectively. Since every pair of edges in $G$ is either adjacent or disjoint, we have $a + b = \binom{k}{2}$. Applying the $2$-edge-balanced condition \eqref{balanced-condition} using the automorphism group orders from Table~\ref{tab:2edge} gives
\[
\frac{b}{a} = \frac{n-3}{4}.
\]
Solving together with $a + b = \binom{k}{2}$ yields
\[
a = \frac{2k(k-1)}{n+1}.
\]
This gives the divisibility condition
\begin{equation}
n + 1 \mid 2k(k-1). \tag{C1}
\end{equation}

We now turn to the constraints imposed by $3$-edge-balance. The five isomorphism classes of graphs in $\mathcal{G}(n,3)$ are listed in Table~\ref{tab:3edge}. Not all of these classes are present for small $n$: $K_{1,3}$ and $P_3\cup K_2$ require $n\ge 5$, and $3K_2$ requires $n\ge 6$. For $n<6$, the absent isomorphism classes do not appear as subgraphs of any graph on $[n]$, so the $3$-edge-balanced condition \eqref{balanced-condition} is imposed only among the classes that are present. For $n=4,5$, a direct check shows that no nontrivial 3-edge-balanced graph exists.

\begin{table}[h]
\centering
\begin{tabular}{llll}
\hline
Graph $H$ & Figure & $\mathrm{Aut}(H)$ & $|\mathrm{Aut}(H)|$ \\
\hline
$K_3$ & \begin{tikzpicture}[baseline=-0.5ex, vertex/.style={circle, fill=black, inner sep=1.5pt}]
\node[vertex] (1) at (0,0) {};
\node[vertex] (2) at (0.5,0) {};
\node[vertex] (3) at (0.25,0.43) {};
\draw (1)--(2)--(3)--(1);
\end{tikzpicture} & $S_3 \times S_{n-3}$ & $6(n-3)!$ \\
$P_4$ & \begin{tikzpicture}[baseline=-0.5ex, vertex/.style={circle, fill=black, inner sep=1.5pt}]
\node[vertex] (1) at (0,0) {};
\node[vertex] (2) at (0.5,0) {};
\node[vertex] (3) at (1,0) {};
\node[vertex] (4) at (1.5,0) {};
\draw (1)--(2)--(3)--(4);
\end{tikzpicture} & $S_2 \times S_{n-4}$ & $2(n-4)!$ \\
$K_{1,3}$ & \begin{tikzpicture}[baseline=-0.5ex, vertex/.style={circle, fill=black, inner sep=1.5pt}]
\node[vertex] (0) at (0.5,0.4) {};
\node[vertex] (1) at (0,0) {};
\node[vertex] (2) at (0.5,0) {};
\node[vertex] (3) at (1,0) {};
\draw (0)--(1);
\draw (0)--(2);
\draw (0)--(3);
\end{tikzpicture} & $S_3 \times S_{n-4}$ & $6(n-4)!$ \\
$P_3 \cup K_2$ & \begin{tikzpicture}[baseline=-0.5ex, vertex/.style={circle, fill=black, inner sep=1.5pt}]
\node[vertex] (1) at (0,0) {};
\node[vertex] (2) at (0.5,0) {};
\node[vertex] (3) at (1,0) {};
\node[vertex] (4) at (1.5,0) {};
\node[vertex] (5) at (2,0) {};
\draw (1)--(2)--(3);
\draw (4)--(5);
\end{tikzpicture} & $S_2 \times S_2 \times S_{n-5}$ & $4(n-5)!$ \\
$3K_2$ & \begin{tikzpicture}[baseline=-0.5ex, vertex/.style={circle, fill=black, inner sep=1.5pt}]
\node[vertex] (1) at (0,0) {};
\node[vertex] (2) at (0.5,0) {};
\node[vertex] (3) at (1,0) {};
\node[vertex] (4) at (1.5,0) {};
\node[vertex] (5) at (2,0) {};
\node[vertex] (6) at (2.5,0) {};
\draw (1)--(2);
\draw (3)--(4);
\draw (5)--(6);
\end{tikzpicture}
 & $(S_2 \wr S_3) \times S_{n-6}$ & $48(n-6)!$ \\
\hline
\end{tabular}
\caption{The five isomorphism classes of graphs in $\mathcal{G}(n,3)$ and their automorphism groups.}
\label{tab:3edge}
\end{table}

Let $c$, $p$, $s$, $d$, $m$ denote the number of subgraphs of $G$ isomorphic to $K_3$, $P_4$, $K_{1,3}$, $P_3 \cup K_2$, $3K_2$, respectively. For $n \ge 6$, applying the $3$-edge-balanced condition \eqref{balanced-condition} using the automorphism group orders from Table~\ref{tab:3edge} gives
\[
 c \cdot |\mathrm{Aut}(K_3)| =  p \cdot |\mathrm{Aut}(P_4)| =  s \cdot |\mathrm{Aut}(K_{1,3})| =  d \cdot |\mathrm{Aut}(P_3\cup K_2)| = m \cdot |\mathrm{Aut}(3K_2)|,
\]
which expresses $p$, $s$, $d$, $m$ in terms of $c$:
\begin{equation}\label{psdm}
\begin{aligned}
p &= 3(n-3)c, \\
s &= (n-3)c, \\
d &= \tfrac{3}{2}(n-3)(n-4)c, \\
m &= \tfrac{1}{8}(n-3)(n-4)(n-5)c.
\end{aligned}
\end{equation}

Now, every $3$-edge subgraph of $G$ is isomorphic to exactly one of the graphs in Table~\ref{tab:3edge}, so
\begin{align*}
\binom{k}{3}
&=c+p+s+d+m \\
&=c\, \left(1+3(n-3)+(n-3)+\frac{3}{2}(n-3)(n-4) +\frac{1}{8}(n-3)(n-4)(n-5)\right) \\
&=c\cdot \frac{(n+1)(n^2-n-4)}{8},
\end{align*}
and hence
\begin{equation}\label{cformula}
c=\frac{4k(k-1)(k-2)}{3(n+1)(n^2-n-4)}.
\end{equation}
This gives the second divisibility condition
\begin{equation}
3(n+1)(n^2-n-4) \mid 4k(k-1)(k-2). \tag{C2}
\end{equation}

We searched over all $n \ge 4$ and all $k$ with $4 \le k \le N/2$ for parameter sets $(n,k)$ that satisfy the divisibility conditions (C1) and (C2) simultaneously, treating $n< 6$ separately as described above. Table~\ref{tab:smallest} gives the 20 smallest parameter sets $(n,k)$ from this search.

\begin{table}[h]
\centering
\begin{tabular}{llll}
\hline
$(13,21)$ & $(21,66)$ & $(34,260)$ & $(77,1378)$ \\
$(84,1275)$ & $(114, 3151)$ & $(139, 3570)$ & $(203, 5644)$ \\
$(229, 11546)$ & $(255, 8321)$ & $(292, 12599)$ & $(309, 2976)$ \\
$(311, 10621)$ & $(414, 37101)$ & $(444, 23941)$ & $(451, 15369)$ \\
$(459, 30706)$ & $(531, 29659)$ & $(539, 49896)$ & $(571, 61776)$ \\
\hline
\end{tabular}
\caption{The 20 smallest admissible parameter sets $(n,k)$ for $3$-edge-balanced graphs.}
\label{tab:smallest}
\end{table}

\section{Search Method and Results}\label{sec4}

We describe in this section our search for $3$-edge-balanced graphs.

For any graph $G \in \mathcal{G}(n,k)$, define its \emph{subgraph count profile} to be the tuple $(c, p, s, d, m)$ counting the number of subgraphs of $G$ isomorphic to $K_3$, $P_4$, $K_{1,3}$, $P_3 \cup K_2$, $3K_2$, respectively. By Lemma~\ref{alltop}, $G$ is $3$-edge-balanced if and only if its subgraph count profile satisfies \eqref{balanced-condition}, which as shown in Section~\ref{sec3}, forces the profile to take a unique value, determined by \eqref{cformula} and \eqref{psdm}. We call this uniquely determined tuple the \emph{target profile} for the parameter pair $(n,k)$, and denote it by $(c^*, p^*, s^*, d^*, m^*)$. Thus $G \in \mathcal{G}(n,k)$ is $3$-edge-balanced if and only if its subgraph count profile is the target profile. The search problem therefore reduces to finding a $k$-edge subgraph of $K_n$ whose subgraph count profile equals $(c^*, p^*, s^*, d^*, m^*)$.

We formulate this as an optimization problem. For $G \in \mathcal{G}(n,k)$ with subgraph count profile $(c, p, s, d, m)$, define its score to be the $\ell^1$ distance from the target profile:
\[
\mathrm{score}(G) = |c - c^*| + |p - p^*| + |s - s^*| + |d - d^*| + |m - m^*|.
\]
Then $G$ is $3$-edge-balanced if and only if $\mathrm{score}(G) = 0$. We search for a zero-score graph using simulated annealing.

The search starts from a near-regular random graph in $\mathcal{G}(n,k)$, constructed by first assigning a target degree of $\lfloor 2k/n \rfloor$ or $\lceil 2k/n \rceil$ to each vertex, then greedily adding edges to meet these requirements, in a random order. Each step picks an edge in $G$ and an edge in $K_n - G$ uniformly at random, and swaps them. The swap is accepted if it does not increase the score, or with probability $\exp(-\Delta / T)$ if it increases the score by $\Delta$, where $T$ is the temperature. The temperature decreases geometrically by a factor $\alpha < 1$ at each step. The search runs for a fixed number of steps, and is restarted from a fresh near-regular random graph in $\mathcal{G}(n,k)$ a fixed number of times.

When a zero-score graph is found, its subgraph count profile is recomputed from scratch to confirm, providing an independent verification that the graph is $3$-edge-balanced.

The search was implemented in C\texttt{++}, with the code developed with the assistance of ChatGPT (OpenAI) \cite{OpenAI:2026}. The program found zero-score graphs for 11 parameter sets $(n,k)=(13,21)$, $(21,66)$, $(34,260)$, $(77,1378)$, $(84,1275)$, $(114,3151)$, $(139,3570)$, $(203,5644)$, $(229,11546)$, $(255,8321)$, $(309,2976)$, including the ten smallest ones. This establishes the existence of nontrivial $3$-edge-balanced graphs. The source code is available at \texttt{https://github.com/ymchee66/t-edge-balanced/}, as are the $3$-edge-balanced graphs found. Table~\ref{tab:statistics} below gives the statistics for the $3$-edge-balanced graphs found.

\begin{table}[ht]
\centering
\begin{tabular}{rrcr}
\hline
\multicolumn{4}{c}{$3$-edge-balanced graph $G$} \\
$n$ & $k$ & $|\mathrm{Aut}(G)|$ & $\lambda$ \\
\hline
13  & 21   & 2 & $15\cdot 10!$ \\
21  & 66   & 1 & $240\cdot 18!$ \\
34  & 260  & 1 & $3552\cdot 31!$ \\
77  & 1378 & 1 & $45792\cdot 74!$ \\
84  & 1275 & 1 & $27930\cdot 81!$ \\
114 & 3151 & 1 & $168840\cdot 111!$ \\
139 & 3570 & 1 & $135456\cdot 136!$ \\
203 & 5644 & 1 & $171864\cdot 200!$ \\
229 & 11546 & 1 & $1025196\cdot 226!$ \\
255 & 8321 & 1 & $277890\cdot 252!$ \\
309 & 2976 & 1 & $7140\cdot 306!$ \\
\hline
\end{tabular}
\caption{The 11 $3$-edge-balanced graphs found, with their parameter sets, automorphism group sizes, and indices $\lambda$.}
\label{tab:statistics}
\end{table}

These examples of $3$-edge-balanced graphs also imply the existence of graphical $3$-$(\binom{n}{2},k,\lambda)$ designs for those $(n,k,\lambda)$ in Table~\ref{tab:statistics}.

\section{The Nonexistence of $t$-Edge-Balanced Graphs for $t\ge 4$}\label{sec5}

We show that no nontrivial $4$-edge-balanced graph exists, which by Lemma~\ref{derived} implies that there no nontrivial $t$-edge-balanced graphs for all $t \ge 4$.

Let $G$ be a $4$-edge-balanced graph of order $n$ and size $k$. By Lemma~\ref{derived}, $G$ is also $3$-edge-balanced, so every isomorphism class of $3$-edge graphs appears in $G$. In particular, we have $c > 0$ and $k > 3$. We now compute $f := n_{G:4K_2}$ in two different ways and derive a contradiction.

For the first computation, we count incidence pairs $(\Delta, F)$ where $\Delta$ is a subgraph of $G$ isomorphic to $K_3$, $F$ is a $4$-edge subgraph of $G$, and $\Delta \subseteq F$. Each triangle in $G$ extends by any of the $k-3$ remaining edges of $G$, so there are $c(k-3)$ such pairs. Counting instead by $F$, the only $4$-edge graphs containing a triangle are $K_3 + e$ (a triangle with one pendant edge) and $K_3 \cup K_2$, each containing exactly one triangle. Hence, letting $x := n_{G:K_3+e}$, we have
\begin{equation}\label{tri-inc}
c(k-3) = x + n_{G:K_3 \cup K_2}.
\end{equation}
Since $G$ is $4$-edge-balanced, Lemma~\ref{alltop} gives $n_{G:H} \cdot |\mathrm{Aut}(H)|$ the same value for all $4$-edge graphs $H$. Therefore
\begin{align*}
n_{G:K_3 \cup K_2} &= x \cdot \frac{|\mathrm{Aut}(K_3+e)|}{|\mathrm{Aut}(K_3 \cup K_2)|} \\
&= x \cdot \frac{2(n-4)!}{12(n-5)!} \\
&= \frac{n-4}{6}\,x.
\end{align*}
Substituting this into \eqref{tri-inc}, we obtain
\begin{equation*}
c(k-3) = \frac{n+2}{6}\,x,
\end{equation*}
or
\begin{equation}\label{x}
x = \frac{6c(k-3)}{n+2}.
\end{equation}
Applying Lemma~\ref{alltop} to $K_3+e$ and $4K_2$, with $|\mathrm{Aut}(4K_2)| = 384(n-8)!$,
\begin{align*}
f &= x \cdot \frac{|\mathrm{Aut}(K_3+e)|}{|\mathrm{Aut}(4K_2)|} \\
&= x \cdot \frac{2(n-4)!}{384(n-8)!} \\
&= \frac{(n-4)(n-5)(n-6)(n-7)}{192}\,x.
\end{align*}
Substituting $x = \frac{6c(k-3)}{n+2}$ gives
\begin{equation}\label{f1}
f = \frac{(n-4)(n-5)(n-6)(n-7)}{32(n+2)}\,c(k-3).
\end{equation}

For the second computation, we count incidence pairs $(M, F)$ where $M$ is a subgraph of $G$ isomorphic to $3K_2$, $F$ is a $4$-edge subgraph of $G$, and $M \subseteq F$. Each copy of $3K_2$ in $G$ extends by any of the $k-3$ remaining edges, giving $m(k-3)$ such pairs. Counting by $F$, the only $4$-edge graphs containing $3K_2$ are $4K_2$, which contains $\binom{4}{3} = 4$ copies of $3K_2$, and $P_3 \cup 2K_2$, which contains two copies of $3K_2$, one for each edge of the path that can be excluded. Hence, letting $y := n_{G:P_3 \cup 2K_2}$, we have
\begin{equation}\label{m(k-3)}
m(k-3) = 4f+2y.
\end{equation}
Applying Lemma~\ref{alltop} to $4K_2$ and $P_3 \cup 2K_2$, with $|\mathrm{Aut}(P_3 \cup 2K_2)| = 16(n-7)!$,
\begin{align*}
y&=f \cdot \frac{|\mathrm{Aut}(4K_2)|}{|\mathrm{Aut}(P_3\cup 2K_2)|} \\
&= f \cdot \frac{384(n-8)!}{16(n-7)!} \\
&= \frac{24}{n-7}\, f.
\end{align*}
Substituting into \eqref{m(k-3)},
\begin{align*}
m(k-3) &= 4f + \frac{48}{n-7}\,f \\
&= \frac{4(n+5)}{n-7}\,f,
\end{align*}
so
\[
f = \frac{n-7}{4(n+5)} m(k-3).
\]

Since $G$ is $3$-edge-balanced, \eqref{psdm} gives $m = \frac{(n-3)(n-4)(n-5)}{8}\,c$, so
\begin{equation}\label{f2}
f = \frac{(n-3)(n-4)(n-5)(n-7)}{32(n+5)}\,c(k-3).
\end{equation}

Since $c(k-3) > 0$, equating \eqref{f1} and \eqref{f2} and cancelling $c(k-3)$ gives
\[
\frac{(n-4)(n-5)(n-6)(n-7)}{n+2} = \frac{(n-3)(n-4)(n-5)(n-7)}{n+5}.
\]
Since $G$ contains $4K_2$ as a subgraph, $n \ge 8$ and we may cancel $(n-4)(n-5)(n-7)$ to obtain
\[
\frac{n-6}{n+2} = \frac{n-3}{n+5},
\]
which is impossible. Hence, no $4$-edge-balanced graphs exist.

\begin{corollary}
There are no nontrivial $t$-edge-balanced graphs for any $t \ge 4$. Equivalently, there are no nontrivial single-orbit graphical $t$-designs for any $t\ge 4$.
\end{corollary}

\section{Conclusion}\label{sec6}

This paper has settled the existence question for $t$-edge-balanced graphs in two directions. For $t = 3$, we have shown that $3$-edge-balanced graphs exist, and found 11 examples of these, including the 10 smallest ones. For $t \ge 4$, we have shown that no nontrivial $t$-edge-balanced graph exists.

Several natural questions remain open. First, are there infinitely many $3$-edge-balanced graphs? In this direction, we make the following conjecture.

\begin{conjecture}
For every pair $(n, k)$ satisfying conditions $\mathrm{(C1)}$ and $\mathrm{(C2)}$, there exists a $3$-edge-balanced graph in $\mathcal{G}(n,k)$.
\end{conjecture}

Second, the 11 examples found here were discovered by computer search. No structural or algebraic explanation for their existence is known. It would be interesting to find explicit constructions, analogous to known construction for $2$-edge-balanced graphs, that would yield infinite families of $3$-edge-balanced graphs.

Finally, by Lemma~\ref{derived}, every $3$-edge-balanced graph is also $2$-edge-balanced. Our 11 examples of $3$-edge-balanced graphs therefore yield new $2$-edge-balanced graphs, since every previously known $2$-edge-balanced graph is triangle-free \cite{Alltop:1966, CaliskanChee:2014, Caliskan:2016}, while our examples contain triangles.

\section*{Acknowledgments}
Research supported by SUTD Grant SKI 2021\_07\_04.

\bibliographystyle{plain}
\bibliography{t-edge-balanced}

\nocite{*}

\end{document}